\documentclass[12pt,leqno]{amsart}
\usepackage{amssymb}
\usepackage{amsmath,amssymb,color}
\numberwithin{equation}{section}

\newcommand{\ep}{\varepsilon}
\newcommand{\ome}{\omega}
\newcommand{\la}{\lambda}
\newcommand{\va}{\varphi}
\newcommand{\wva}{\widehat{\varphi}}
\newcommand{\ppp}{\partial}
\newcommand{\wh}{\widehat}

\newcommand{\ddd}{\mbox{div}\thinspace}
\newcommand{\rrr}{\mbox{rot}\thinspace}

\newcommand{\weight}{e^{2s\alpha}}

\newcommand{\R}{\mathbb{R}}

\newcommand{\N}{\mathbb{N}}

\newcommand{\www}{\widetilde}

\newcommand{\ooo}{\overline}
\newcommand{\OOO}{\Omega}

\allowdisplaybreaks
\title
[]
{Global Lipschitz stability for an inverse source problem for
the Navier-Stokes equations}
\author{
$^1$ Oleg Y. Imanuvilov, and $^{2,3,4,5}$ M.~Yamamoto }

\thanks{
$^1$ Department of Mathematics, Colorado State University\\
101 Weber Building, Fort Collins, CO 80523-1874, USA\\
e-mail: {\tt oleg@math.colostate.edu}\\
$^2$ Graduate School of Mathematical Sciences, The University
of Tokyo, Komaba, Meguro, Tokyo 153-8914, Japan \\
$^3$ Honorary Member of Academy of Romanian Scientists, 
Ilfov, nr. 3, Bucuresti, Romania \\
$^4$ Correspondence member of Accademia Peloritana dei Pericolanti,\\
Palazzo Universit\`a, Piazza S. Pugliatti 1 98122 Messina Italy \\
$^5$ Peoples' Friendship University of Russia 
(RUDN University) 6 Miklukho-Maklaya St, Moscow, 117198, Russian Federation
e-mail: {\tt myama@ms.u-tokyo.ac.jp}
}

\date{}
\begin{document}
\maketitle
\begin{abstract}
For linearized Navier-Stokes equations, we consider an inverse source 
problem of determining a spatially varying divergence-free factor.
We prove the global Lipschitz stability by interior data over a time 
interval and velocity field at $t_0>0$ over the spatial domain.
The key are Carleman estimates for the Navier-Stokes equations and 
the operator rot.
\end{abstract}

\section{Introduction}

Let $T > 0$, and $\OOO \subset \R^3$ be a bounded domain with smooth boundary
$\ppp\OOO$ and $A$, $B$ be sufficiently smooth as specified below.
By $a^T$ we denote the transpose of a vector $a$ under consideration.
We set 
$$
\ppp_j = \frac{\ppp}{\ppp x_j}, \quad j=1,2,3, \quad 
\ppp_t = \frac{\ppp}{\ppp t}, \quad 
\nabla = (\ppp_1, \ppp_2, \ppp_3), \quad \nabla_{x,t} = (\nabla, \ppp_t)
$$
and
$$
\ddd v = \sum_{j=1}^3 \ppp_jv_j \quad \mbox{for $v=(v_1,v_2,v_3)^T$}.
$$
 
%For a linearized Navier-Stokes equations, we establish a Carleman 
%estimate only by conventional $L^2$-Carleman estimates for a parabolc and
%the Laplace equation.
We consider
$$
\left\{\begin{array}{rl}
& \ppp_tv(x,t) - \Delta v + (A(x,t)\cdot\nabla)v
+ (v\cdot \nabla)B(x,t) + \nabla p = F(x,t),\\        
& \ddd v = 0 \quad \mbox{in $\OOO\times (0,T)$}, \quad v=0
\quad \mbox{on $\ppp\OOO \times (0,T)$}. 
\end{array}\right.
                                                         \eqno{(1.1)}
$$

In this article, we discuss the following inverse source problem when a 
subdomain $\omega \subset \OOO$ and $t_0 \in (0,T)$ are given.
\\
\vspace{0.2cm}
{\bf Inverse source problem}.
{\it 
We assume that a source term $F$ satisfies some conditions stated below.  
Then determine $F$ by $v$ in $\ome \times (0,T)$ and $v(\cdot,t_0)$
in $\OOO$.
}
\\

In particular, we aim at the global Lipschitz stability for the inverse 
problem.

We notice that in our inverse source problem, we do not assume any 
data of the pressure field $p$, because the observation of $p$ may 
be more difficult than the velocity field $v$.
If suitable data of $p$ are available, then the inverse problem 
is easily solved in the same way as Imanuvilov and Yamamoto 
\cite{IY1}, and we omit the details.

We introduce some notations.  We set 
$$
Q:= \OOO \times (0,T).
$$
%Let $Q\subset \R^4$ be a domain in $\{(x,t);\, x\in \R^3, \, t>0\}$,
Let $\gamma = (\gamma_1, \gamma_2, \gamma_3) \in (\N \setminus \{0\})^3$ 
and $\vert \gamma\vert = \gamma_1 + \gamma_2 + \gamma_3$, and 
$\ppp_x^{\gamma} = \ppp_1^{\gamma_1}\ppp_2^{\gamma_2}\ppp_3^{\gamma_3}
\\
= \left( \frac{\ppp}{\ppp x_1}\right)^{\gamma_1} 
\left(\frac{\ppp}{\ppp x_2}\right)^{\gamma_2} 
\left(\frac{\ppp}{\ppp x_3}\right)^{\gamma_3}$.
We further set
$$
H^{2,1}(Q) = \{ v \in L^2(Q); \thinspace
\ppp_x^{\gamma}v \in L^2(Q), \thinspace \vert \gamma\vert \le 2,
\thinspace \ppp_tv \in L^2(Q)\},
$$
$$
H:= \ooo{ \{ v\in C^{\infty}_0(\OOO)^3:\, \ddd v = 0 \, \quad \mbox{in $\OOO$}
\}^{L^2(\OOO)} }
$$
and
$$
V:= (H^1_0(\OOO))^3 \cap H.
$$

In this article, for simplicity, we assume some smoothness for 
$v, p$, and not pursue the minimum regularity.
More precisely, we assume
$$
A \in C^2([0,T];H^1(\OOO)), \quad B\in C^2([0,T];H^2(\OOO)).
$$
For the solution $v, p$, we assume the regularity:
$$
\ppp_t^kv \in L^2(0,T;V) \cap H^{2,1}(Q), \quad k=0,1,2, \quad
p \in L^2(0,T;H^1(\OOO)).               \eqno{(1.2)}
$$
\\

Let $\ome \subset \OOO$ be an arbitrarily chosen non-empty subdomain.
For an arbitrarily fixed constant $M > 0$, we define an admissible set of 
unknown $F$ by  
$$
\mathcal{F} := \{ F \in H^2(0,T;H^2(\OOO));\, \ddd F(x,t_0) = 0 \,\,
\mbox{in $\OOO$},
$$
$$
\vert \ppp_t^k F(x,t)\vert \le M\vert F(x,t_0)\vert \,\,\mbox{for 
$(x,t) \in Q$ and $k=1,2$},\quad 
F\vert_{\ome} = F\vert_{\ppp\OOO} = 0.\}        \eqno{(1.3)}
$$
\\

We are ready to state our main result.
\\
{\bf Theorem 1.}\\
{\it
There exists a constant $C>0$, depending on $\mathcal{F}$, such that 
$$
\Vert F\Vert_{H^2(0,T;L^2(\OOO))} 
\le C(\Vert v\Vert_{H^2(0,T;H^1(\ome))} + \Vert v(\cdot,t_0)\Vert_{H^2(\OOO)}
)
$$
for each $F \in \mathcal{F}$.
}
\\

Theorem 1 asserts the global Lipschitz stability over $\OOO$ for the
inverse problem.  Theorem 1 immediately implies the uniqueness within 
$\mathcal{F}$.
In general we can prove the uniqueness within more generous 
admissible set than $F$ defined by (1.3).  For better uniqueness,
we need different arguments, especially another type of Carleman estimate,
and we omit the details.

In Choulli, Imanuvilov, Puel and Yamamoto \cite{CIPY}, a similar 
stability estimate for the same type of the inverse source pronlem was proved.
Theorem 1 provides the global Lipschitz stability under conditions on 
$R$ and $\omega$ which are more natural and easier to be verified. 
\\

%{\bf Remark 1.}\\

The condition 
$$
\ddd F(x,t_0) = 0 \quad \mbox{in $\OOO$}               \eqno{(1.4)}
$$
is essential for the uniqueness in determining $F$ without data of $p$, 
as the following example 
shows.
\\
{\bf Obstruction to the uniqueness:}
\\
We consider a simple case
$$
\left\{\begin{array}{rl}
& \ppp_tv(x,t) - \Delta v + \nabla p = f(x),  \\
& \ddd v = 0, \quad (x,t) \in Q:= \OOO \times I,\\
& v(x,t_0) = 0 \quad\mbox{for $x\in \OOO$}, \quad
\mbox{supp $v \subset \OOO \times (0,T)$}. 
\end{array}\right.
                                                         \eqno{(1.5)}
$$
Here $f$ is an $\R^3$-valued smooth function.
It is trivial that $(v,p) = (0,0)$ satisfies (1.5) with $f=0$.
Let $\psi\in C^{\infty}_0(\OOO)$.  Then $(v,p) := (0, \, \psi)$
satisfies (1.5) with $f:= \nabla \psi$.  In other words, by the 
appearance of the pressure field $p$ in the Navier-Stokes equations, 
there is no possibility for the uniqueness in 
determining $f$ given by a scalar potential.
\\

We notice that if $F(x,t_0)$ is given by a vector potential:
$$
F(x,t_0) = \rrr q(x), \quad x \in \OOO
$$
with smooth $q$, then (1.4) holds automatically.
\\

Next we examine the condition
$$
\vert \ppp_t^k F(x,t)\vert \le M\vert F(x,t_0)\vert \quad\mbox{for 
$(x,t) \in Q$ and $k=1,2$}.                               \eqno{(1.6)}
$$
{\bf Examples}.\\
{\bf (i)} Let $R(x,t)$ be a smooth $3\times 3$ matrix such that 
det $R(x,t_0) \ne 0$ for $x \in\ooo{\OOO}$.  Then (1.6) is satisfied.
\\
{\bf (ii)}  Let $r(x,t) = (r_1(x,t), r_2(x,t), r_3(x,t))^T$ and 
$f(x) = (f_1(x),\, f_2(x), \, 0)^T$ be smooth functions satisfying 
$$
\rrr f(x) = (-\ppp_3f_2(x), \, \ppp_3f_1(x),\, \ppp_1f_2(x) - \ppp_2f_1(x))
^T = 0, \quad
\rrr r(\cdot,t_0) = 0 \quad \mbox{in $\OOO$}.
$$
Then $F(x,t):= r(x,t) \times f(x)$ satisfies (1.4).  
Indeed, as is directly seen, 
$$
\ddd F(x,t_0) = \ddd (r(x,t_0) \times f(x)) 
= f(x)\cdot \rrr r(x,t_0) - r(x,t_0)\cdot \rrr f(x) = 0 \quad
\mbox{in $\OOO$.}
$$
If we assume $r_3(x,t_0) \ne 0$ for $x \in \ooo{\OOO}$, then (1.6) is 
satisfied.  indeed, by 
$$
r \times f = (-r_3f_2, \, r_3f_1, \, r_1f_2 - r_2f_1)^T \quad 
\mbox{in $Q$},
$$
we have
$$
\vert \ppp_t^k(r\times f)(x,t)\vert \le C(\vert f_1(x)\vert 
+ \vert f_2(x)\vert), \quad (x,t)\in Q, \, k=1,2
$$
and
$$
\vert r(x,t_0) \times f(x)\vert 
\ge \vert r_3(x,t_0)\vert (\vert f_1(x)\vert + \vert f_2(x)\vert)
\ge C(\vert f_1(x)\vert + \vert f_2(x)\vert), \quad x\in \OOO
$$
by $r_3(x,t_0) \ne 0$ for $x \in \ooo{\OOO}$.
\\

Our proof is based on Bukhgeim and Klibanov \cite{BK} for single 
partial differential equations.  The work relies essentially on 
an $L^2$-weighted estimate called a Carleman estimate.
For the Navier-Stokes equations, we need such a Carleman estimate 
and we here modify the Carleman estimate 
proved in Choulli, Imanuvilov, Puel and Yamamoto \cite{CIPY}.

As for inverse problems by Carleman estimates after \cite{BK},
we refer to very limited works: Beilina and Klibanov \cite{BeKl},
Bellassoued and Yamamoto \cite{BY}, Huang, Imanuvilov and 
Yamamoto \cite{HIY}, Imanuvilov and Yamamoto \cite{IY1}, 
\cite{IY2}, Isakov \cite{Is}, Klibanov \cite{Kl},
Klibanov and Timonov \cite{KT}, Yamamoto \cite{Ya}.
Epecially for related works on inverse problems for the incompressible
fluid equations,
see Bellassoued, Imanuvilov and Yamamoto \cite{BIY}, 
Choulli, Imanuvilov, Puel and Yamamoto \cite{CIPY}, 
Fan, Di Cristo, Jiang and Nakamura \cite{FDJN},
Fan, Jiang and Nakamura \cite{FJN}, 
Imanuvilov, Lorenzi and Yamamoto \cite{ILY}, Imanuvilov and 
Yamamoto \cite{IYNS}.
Except for \cite{BIY} and \cite{ILY}, all these articles use the same type of 
Carleman estimate whose weight has singularities at both ends of the time
interval.
See Imanuvilvov and Yamamoto \cite{IY3} for compressible fluid equations, and
also Fern\'andez-Cara, Guerrero, Imanuvilov
and Puel \cite{FGIP} for some exact controllability problem for the 
Navier-Stokes equations.

The article is composed of five sections.  In Section 2, as Lemmata 
1 and 2, we show two key Carleman estimates for the Navier-Stokes equations
and for stationary Maxwell's equations.
In Section 3, we complete the proof of Theorem 1.
Sections 4 and 5 are devoted respectively to concluding remarks and 
the derivation of the first key Carleman estimate with general powers.

\section{Key Carleman estimates}

For a non-empty suddomain $\ome_0$ with $\ooo{\ome_0} \subset \ome$,
we know (e.g., Imanuvilov \cite{Ima}) that there exists $\eta \in
C^2(\ooo{\OOO})$ satisfying
$$
\eta\vert_{\ppp\OOO} = 0, \quad \eta > 0 \quad \mbox{in $\OOO$}, \quad 
\vert \nabla \eta \vert > 0 \quad \mbox{on $\ooo{\OOO\setminus \ome_0}$}.
                                             \eqno{(2.1)}
$$

Henceforth we set
$$
t_0 = \frac{T}{2}.
$$
Let $\ell \in C^{\infty}[0,T]$ satisfy
$\ell(t) = \ell(T-t)$ for $0\le t \le T$, and
$$
\left\{ \begin{array}{rl}
& \ell(t) > 0, \quad 0<t<T, \\
& \ell(t) = t, \quad 0\le t\le \frac{T}{4}, \\
&\ell(t_0) > \ell(t), \quad t \in (0,T) \setminus \{ t_0\}.
\end{array}\right.
                                               \eqno{(2.2)}
$$
We set 
$$
Q_{\ome}: = \ome \times (0,T),
$$
recall that $Q = \OOO \times (0,T)$, and
$$
\va(x,t) := \frac{e^{\la\eta(x)}}{\ell^8(t)}, \quad
\alpha(x,t) := \frac{e^{\la\eta(x)} - e^{2\la\Vert\eta\Vert_{C(\ooo{\OOO})}}}
{\ell^8(t)}, \quad (x,t) \in Q.                 \eqno{(2.3)}
$$
Here we fix a sufficiently large constant $\la>0$.

Hencforth $C>0$, $C_j> $ denote generic constants which are dependent on 
$\OOO$, $T$, $\la$, but independent of $s$.
A Carleman estimate for the Navier-Stokes equations is stated as:
\\
{\bf Lemma 1.}\\
{\it
Let $F \in L^2(0,T;H)$, $v(\cdot,t) \in V$, $0\le t\le T$ and 
$v \in L^2(0,T;V) \cap H^{2,1}(Q)$ satisfy (1.1).  Let 
$m\in \N \cup \{0\}$.  Then there exist constants 
$C>0$ and $s_0>0$ such that 
\begin{align*}
& \int_Q (s^m\va^m\vert \nabla v\vert^2 + s^{m+1}\va^{m+1}\vert \rrr v\vert^2
+ s^{m+2}\va^{m+2}\vert v\vert^2) \weight dxdt \\
\le & C\int_Q s^m\va^m\vert F\vert^2 \weight dxdt\\
+ & C\int_{Q_{\ome}} (s^{m+1}\va^{m+1}\vert \rrr v\vert^2 
+ s^{m+2}\va^{m+2}\vert v\vert^2
+ s^{m+1}\va^{m+1}\vert \nabla v\vert^2) \weight dxdt 
\end{align*}
for all $s\ge s_0$.
}
\\

For the case of $m=0$, the proof is found e.g., Theorem 2 in 
Choulli, Imanuvilov, Puel and Yamamoto \cite{CIPY}.  
In Section 5, for completeness we derive the lemma with general $m\in \N$
from the case $m=0$.
\\

Next let $\psi \in C^2(\ooo{\OOO})$ satisfy 
$$
\psi > c_0 \quad \mbox{in $\OOO$}, \quad \psi = c_0 \quad 
\mbox{on $\ppp\OOO$},\quad \vert \nabla \psi\vert > 0 \quad
\mbox{on $\ooo{\OOO\setminus \ome}$},                        \eqno{(2.4)}
$$
where $c_0$ is a given constant.
With sufficiently large fixed $\la > 0$, we set 
$$
\va_0(x) := e^{\la\psi(x)}, \quad x\in \OOO.
$$
Then we show
\\
{\bf Lemma 2.}\\
{\it
Let $w \in H^2(\OOO) \cap H^1_0(\OOO)$ satisfy $w\vert_{\ome} = 0$, and
$$
\rrr w = g, \quad \ddd w = h \quad \mbox{in $\OOO$}.   \eqno{(2.5)}
$$
Then there exist constants $C>0$ and $s_1>0$ such that 
$$
\int_\OOO \left( \frac{1}{s}\vert \nabla w\vert^2 + s\vert w\vert^2
\right) e^{2s\va_0(x)} dx 
\le C\int_\OOO (\vert g\vert^2 + \vert h\vert^2) e^{2s\va_0(x)} dx
$$
for all $s\ge s_1$.
}
\\

This is a Carleman estimate for stationary Maxwell's equations and we
can refer to Vogelsang \cite{Vo}.  Here we prove by means of 
a Carleman estimate for $-\Delta$, not by \cite{Vo}.
\\
{\bf Proof of Lemma 2.}\\
Setting $\www{\psi}(x) = \psi(x) - c_0$, $x\in \OOO$, we see that 
$\www{\psi} > 0$ in $\OOO$, $\psi = 0$ on $\ppp\OOO$ and 
$\vert \nabla \www{\psi}\vert > 0$ on $\ooo{\OOO \setminus \ome}$.
We set $\www{\va_0}(x):= e^{\la\www{\psi}(x)}
= e^{-\la c_0}\va_0(x)$, $x\in \OOO$.

Since 
$$
\rrr \rrr w = -\Delta w + \nabla (\ddd w)
$$
by $w \in H^2(\OOO)$,  we see
$$
-\Delta w = \rrr g - \nabla h \quad \mbox{in $\OOO$}     \eqno{(2.6)}
$$
by (2.5).  In terms of $w\vert_{\ppp\OOO} = 0$ and $w\vert_{\ome} = 0$, we 
apply an $H^{-1}$-Carleman esimate (e.g., Theorem 2.2 in Imanuvilov and 
Puel \cite{IP}) to (2.6), and we obtain
$$
\int_\OOO \left( \frac{1}{s}\vert \nabla w\vert^2 + s\vert w\vert^2
\right) e^{2s\www{\va_0}(x)} dx 
\le C\int_\OOO (\vert g\vert^2 + \vert h\vert^2) e^{2s\www{\va_0}(x)} dx
$$
for all $s\ge \www{s_1}$: some constant.
Setting $\www{s}:= se^{\la c_0}$, we have
$$
\int_\OOO \left( \frac{1}{\www{s}}\vert \nabla w\vert^2 + \www{s}\vert w\vert^2
\right) e^{2\www{s}\www{\va_0}(x)} dx 
\le C\int_\OOO (\vert g\vert^2 + \vert h\vert^2) e^{2\www{s}\www{\va_0}(x)} dx
$$
for $\www{s} \ge \www{s_1}$.
By $\www{s}\www{\va_0}(x) = se^{\la c_0}e^{-\la c_0}\va_0(x)
= s\va_0(x)$ for $x\in\OOO$ and $se^{\la c_0} \ge s$, we see
$$
\int_\OOO \left( \frac{e^{-\la c_0}}{s}\vert \nabla w\vert^2 + s\vert w\vert^2
\right) e^{2s\va_0(x)} dx 
\le C\int_\OOO (\vert g\vert^2 + \vert h\vert^2) e^{2s\va_0(x)} dx
$$
for $s \ge \www{s_1}e^{-\la c_0}$.  Setting $s_1 := \www{s_1}e^{-\la c_0}$
and replacing $C$ by $Ce^{\la c_0}$, we complete the proof of 
Lemma 2.   $\blacksquare$
\\

We conclude this section with a lemma which constitutes a technical part of 
the proof of Theorem 1, but the proof is elementary.
\\
{\bf Lemma 3.}\\
{\it
Let $g \in L^1(\OOO)$.  Then
$$
\int_Q \va(x,t) \vert g(x)\vert e^{2s\alpha(x,t)} dxdt 
\le C\int_{\OOO} \vert g(x)\vert e^{2s\alpha(x,t_0)} dx,
$$
where $\alpha$ and $\va$ are defined by (2.3).
}
\\
{\bf Proof of Lemma 3.}
\\
We set 
$$
\wh{\va}(t) = \frac{1}{\ell^8(t)}, \quad 0<t<T.    \eqno{(2.7)}
$$
Then, as can be directly verified, fixing 
$\la>0$, we can choose a constant $C>0$ such that 
$$
C^{-1} \va(x,t) \le \wh{\va}(t) \le C\va(x,t), \quad (x,t) \in Q.
                                                      \eqno{(2.8)}
$$
By (2.8), it is sufficient to prove 
$$
\int_Q \wh{\va}(t) \vert g(x)\vert e^{2s\alpha(x,t)} dxdt 
\le C\int_{\OOO} \vert g(x)\vert e^{2s\alpha(x,t_0)} dx.
$$
First 
$$
\int_Q \wh{\va}(t) \vert g(x)\vert e^{2s\alpha(x,t)} dxdt 
= \int_{\OOO} \vert g(x)\vert e^{2s\alpha(x,t_0)}
\left( \int^T_0 \wh{\va}(t) e^{2s(\alpha(x,t)-\alpha(x,t_0))} dt\right) dx.
$$
Since $\wh\va(t) \ge 0$ for $0<t<T$, we have
\begin{align*}
& \int^T_0 \wh{\va}(t) e^{2s(\alpha(x,t)-\alpha(x,t_0))} dt
\le \int^T_0 \wh{\va}(t)\exp\left( -2sC\frac{\ell^8(t_0)-\ell^8(t)}
{\ell^8(t_0)\ell^8(t)} \right) dt\\
\le& \int^T_0 \wh{\va}(t)\exp\left( -\mu_1s\frac{\ell^8(t_0)-\ell^8(t)}
{\ell^8(t)} \right) dt,
\end{align*}
with some constant $\mu_1 > 0$.
We fix $\ep>0$.  Hence
\begin{align*}
& \int^T_0 \wh{\va}(t)e^{2s(\alpha(x,t) - \alpha(x,t_0))} dt\\
\le & \left( \int^{\ep}_0 + \int^T_{T-\ep} 
+ \int^{T-\ep}_{\ep}\right)
\frac{1}{\ell^8(t)} \exp\left( -\mu_1s\frac{\ell^8(t_0)-\ell^8(t)}
{\ell^8(t)} \right) dt.
\end{align*}
Here we have
$$
\int^{\ep}_0 \frac{1}{\ell^8(t)} 
\exp\left( -\mu_1s\frac{\ell^8(t_0)-\ell^8(t)}{\ell^8(t)} \right) dt
\le \int^{\ep}_0 \frac{1}{\ell^8(t)} 
e^{-\frac{C}{\ell^8(t)}} dt < \infty
$$
for $s\ge 1$ for example.
Moreover 
$$
\int^{T-\ep}_ {\ep} \frac{1}{\ell^8(t)} 
\exp\left( -\mu_1s\frac{\ell^8(t_0)-\ell^8(t)}{\ell^8(t)} \right) dt
\le \int^{T-\ep}_{\ep} \frac{1}{\ell^8(t)} dt < \infty.
$$
Thus the proof of Lemma 3 is complete.
\section{Proof of Theorem 1}

With loss of generality, we can assume that $t_0 = \frac{T}{2}$.
Indeed, for arbitrary $t_0\in (0,T)$, we choose small $\delta>0$
such that $0<t_0-\delta < t_0+\delta < T$.  Then, in place of $(0,T)$,
we can continue the proof over the time interval $(t_0-\delta, \,
t_0+\delta)$.
We set
$$
D := \Vert v\Vert_{H^2(0,T;H^1(\ome))} + \Vert v(\cdot,t_0)\Vert_{H^2(\OOO)}.
$$
Applying Lemma 1 with $m=1$ to $\ppp_tv$ and $\ppp_t^2v$, we have
\begin{align*}
& \int_Q s^2\va^2(\vert \rrr \ppp_tv\vert^2 + \vert \rrr\ppp_t^2v\vert^2)
\weight dxdt 
\le C\int_Q s\va \sum_{k=1}^2 \vert \ppp_t^kF\vert^2 \weight dxdt\\
+ & CD^2\int_{Q_{\ome}} s^3\va^3 \weight dxdt 
\end{align*}
for all $s\ge s_0$.
Since $\int_Q s^3\va^3 \weight dxdt \le C$ and 
$\sum_{k=1}^2 \vert \ppp_t^kF(x,t)\vert^2 \le M\vert F(x,t_0)\vert$ for 
$(x,t)\in Q$, we obtain
$$
\int_Q s\va^2(\vert \rrr \ppp_tv\vert^2 + \vert \rrr\ppp_t^2v\vert^2)
\weight dxdt 
\le C\int_Q \va \vert F(x,t_0)\vert^2 e^{2s\alpha(x,t)} dxdt
+ CD^2.                                    \eqno{(3.1)}
$$
By (1.1), we have
\begin{align*}
& F(x,t_0) = \ppp_tv(x,t_0) - \Delta v(x,t_0) \\
+ & (A(x,t_0)\cdot\nabla)v(x,t_0) + (v(x,t_0)\cdot\nabla)B(x,t_0) 
+ \nabla p(x,t_0).
\end{align*}
Since $\rrr \nabla = 0$, we obtain
$$
\rrr F(x,t_0) = \rrr \ppp_tv(x,t_0) + \rrr a(x), \quad x\in \OOO,
                                                      \eqno{(3.2)}
$$
where 
$$
a(x):= - \Delta v(x,t_0) + (A(x,t_0)\cdot\nabla)v(x,t_0) 
+ (v(x,t_0)\cdot\nabla)B(x,t_0), \quad x \in \OOO.
$$
Hence, using that $\alpha(x,0) = \lim_{t\downarrow 0} \alpha(x,t) = 0$
for $x\in \OOO$, we have
\begin{align*}
& \int_{\OOO} \vert \rrr \ppp_tv(x,t_0)\vert^2 e^{2s\alpha(x,t_0)} dx
= \int^{t_0}_0 \ppp_t \left( \int_{\OOO}
\vert \rrr \ppp_tv(x,t)\vert^2 e^{2s\alpha(x,t)} dx \right) dt\\
=& 2\int^{t_0}_0\int_{\OOO} ((\rrr\ppp_tv \cdot \rrr\ppp_t^2v)
+ \vert \rrr \ppp_tv(x,t)\vert^2 s(\ppp_t\alpha)) e^{2s\alpha(x,t)} dx dt.
\end{align*}
Henceforth we write $\ell'(t) = \frac{d \ell}{dt}(t)$.
Since
$$
\ppp_t\alpha(x,t) = \frac{-8}{\ell^8(t)}\frac{\ell'(t)}{\ell(t)}
(e^{\la\eta(x)} - e^{\la\Vert\eta\Vert_{C(\ooo{\OOO})}})
$$
and
$$
\left\vert \frac{\ell'(t)}{\ell(t)}\right\vert \le \frac{C}{\ell^8(t)},
\quad 0<t<T,
$$
we can verify that $\vert \ppp_t\alpha(x,t)\vert \le C\va(x,t)^2$ for 
$(x,t) \in Q$.  Hence, applying $\vert ab\vert \le \frac{1}{2}
(\vert a\vert^2 + \vert b\vert^2)$ for $a,b \in \R$, we obtain
\begin{align*}
& \int_{\OOO} \vert \rrr \ppp_tv(x,t_0)\vert^2 e^{2s\alpha(x,t_0)} dx
\le C\int_Q (\vert \rrr\ppp_tv\vert \vert \rrr\ppp_t^2v\vert
+ s\va^2\vert \rrr\ppp_tv\vert^2) e^{2s\alpha(x,t)} dxdt\\
\le& C\int_Q (\vert \rrr\ppp_tv\vert^2 + \vert \rrr\ppp_t^2v\vert^2
+ s\va^2\vert \rrr\ppp_tv\vert^2) e^{2s\alpha(x,t)} dxdt
\end{align*}
for all large $s>0$.  Since $\va > 0$ on $\ooo{Q}$, choosing $s>0$ large,
we can estimate
$$
\vert \rrr\ppp_tv\vert^2 \le Cs\va^2 \vert \rrr\ppp_tv\vert^2
$$
in $Q$.  Therefore, application of (3.1) yields
\begin{align*}
& \int_{\OOO} \vert \rrr \ppp_tv(x,t_0)\vert^2 e^{2s\alpha(x,t_0)} dx
\le C\int_Q (\vert \rrr \ppp_t^2v\vert^2 
+ s\va^2\vert \rrr \ppp_tv(x,t)\vert^2)
e^{2s\alpha(x,t)} dxdt\\
\le & C\int_Q \va(x,t)\vert F(x,t_0)\vert^2 e^{2s\alpha(x,t)} dxdt
+ CD^2
\end{align*}
for all large $s>0$.
Applying Lemma 3 to the first term on the right-hand side, we have
$$
\int_{\OOO} \vert \rrr \ppp_tv(x,t_0)\vert^2 e^{2s\alpha(x,t_0)} dx
\le C\int_{\OOO} \vert F(x,t_0)\vert^2 e^{2s\alpha(x,t_0)} dx
+ CD^2.
$$
Therefore (3.2) yields
\begin{align*}
& \int_{\OOO} \vert \rrr F(x,t_0)\vert^2 e^{2s\alpha(x,t_0)} dx\\
\le &C\int_{\OOO} \vert F(x,t_0)\vert^2 e^{2s\alpha(x,t_0)} dx
+ CD^2
+ C\int_{\OOO} \vert \rrr a(x)\vert^2 e^{2s\alpha(x,t_0)} dx \\
\le & C\int_{\OOO} \vert F(x,t_0)\vert^2 e^{2s\alpha(x,t_0)} dx
+ CD^2
\end{align*}
for all large $s>0$.  Applying Lemma 2 to $F(x,t_0)$ in view of 
$\ddd F(x,t_0) = 0$ for $x\in \OOO$ and $F\vert_{\ppp\OOO} 
= F\vert_{\ome} = 0$, we obtain
$$
\int_{\OOO} s\vert F(x,t_0)\vert^2 e^{2s\alpha(x,t_0)} dx
\le C\int_{\OOO} \vert \rrr F(x,t_0)\vert^2 e^{2s\alpha(x,t_0)} dx.
$$
for all large $s>0$.   Choosing $s>0$ large, we can complete the proof of 
Theorem 1.
$\blacksquare$

\section{Concluding Remarks}

{\bf 4.1. Nonlinear case.}

In this article, we mainly consider linearized Navier-Stokes equations.
For the original Navier-Stokes equations
$$
\ppp_tv - \Delta v + (v\cdot v)v + \nabla p = F, \quad 
\ddd v = 0 \quad \mbox{in $Q$},              \eqno{(4.1)}
$$
assuming suitable smoothness of $v, p$ and taking the difference
$v_1 - v_2$ where $v_1$ and $v_2$ are solutions to (4.1) respectively 
with $F_1$ and $F_2$, we can reduce the inverse source problem to 
the one for the linearized Navier-Stokes equations.  We omit the details.

{\bf 4.2. Case of $t_0=0$.}

We do not assume any initial value but we are given $v(\cdot,t_0)$ with
$t_0 > 0$.  In other words, our inverse problem does nor correspond 
to an initial boundary value problem for the Navier-Stokes equations.
The formulation of inverse problems with $t_0>0$ is standard for 
parabolic equations (e.g., \cite{BIY}, \cite{CIPY},
\cite{FDJN}, \cite{FJN}, \cite{HIY}, 
\cite{ILY}, \cite{IY1}, \cite{Is}, \cite{Ya}).

On the other hand, for the case of $t_0=0$, inverse problems of 
the type discussed here for equations of parabolic type including 
the Navier-Stokes equations, are open, 
and even the uniqueness is unknown in general.
In some special cases, we can prove the uniqueness (e.g., Theorem 4.7 in 
\cite{Kl}).

{\bf 4.3. Arbitrariness of observation time interval.}

We consider the inverse problem in the time interval $(0,T)$, when
the Navier-Stokes equations (1.1) hold. 
As is seen by the proof, in view of the parabolicity of the equations,
it is not necessary to take data over the same time interval, but sufficient
over a smaller time interval $I$ such that $\ooo{I} \subset (0,T)$.
Thus for the inverse problem, it suffices that we are concerned with 
the smoothness of the solution $v, p$ for 
$t>0$, and not to exploit the regularity at $t=0$ which is more 
delicate for $t>0$.

{\bf 4.4. Choices of the weight functions of Carleman estimates.}
 
Our key is a Carleman estimate.  We have two kinds of Carleman estimates
according to the weight functions:
\begin{itemize}
\item
Regular weight function:
$$
\va(x,t) := e^{\la(d(x) - \beta (t-t_0)^2)}, \quad (x,t) \in Q.    \eqno{(4.2)}
$$
\item
Singular weight function:
$$
\va(x,t) = \frac{e^{\la\eta(x)} - e^{2\la\Vert \eta\Vert_{C(\ooo{\OOO})}}}
{h(t)},          \quad (x,t) \in Q,                     \eqno{(4.3)}
$$
where $\lim_{t\to 0} h(t) = \lim_{t\to T} h(t) = 0$. 
\end{itemize}
Here $d, \eta \in C^2(\ooo{\OOO})$ are chosen suitably.

The Carleman estimate with the weight (4.3) was proved firstly by 
Imanuvilov \cite{Ima} for a parabolic equation.
As for related inverse problems for the Navier-Stokes equations,
Choulli, Imanuvilov, Puel and Yamamoto \cite{CIPY},
Fan, Di Cristo, Jiang and Nakamura \cite{FDJN}, 
Fan, Jiang and Nakamura \cite{FJN} used Carleman estimates with 
(4.3), while Bellassoued, Imanuvilov and Yamamoto \cite{BIY},
Imanuvilov, Lorenzi and Yamamoto \cite{ILY} and the 
current article rely on (4.2).

The Carleman estimate with (4.2) is relevant for considering  
inverse problems locally in $x$ in the case where we are not given
the boundary condition on the whole $\ppp\OOO$, while the Carleman 
estimate with (4.3) is more convenient for deriving the global
stability in inverse problems.

\section{Appendix. Proof of Lemma 1}

We derive the lemma from the case of $m=0$ which is proved in 
\cite{CIPY}.  We reacall that $\wva$ is defined by (2.7).
We set 
$$
w(x,t) := \wva^{\frac{m}{2}}(t)v(x,t), \quad (x,t) \in Q.
$$
Then 
$$
\ppp_tw = \wva^{\frac{m}{2}}\ppp_tv + \frac{m}{2}\wva^{\frac{m}{2}-1}
\wva'(t)v(x,t)
= \wva^{\frac{m}{2}}\ppp_tv + q(t)\wva w
$$
with some $q\in L^{\infty}(0,T)$, by noting 
$$
\wva^{\frac{m}{2}-1}\wva'(t)v = \wva^{-1}\wva'(t)w 
= \wva(\wva^{-2}\wva'(t))w = -8\ell'(t)\ell^7\wva w
$$
and $\ell'\ell^7 \in L^{\infty}(0,T)$.

For simplicity, we set 
$$
L(v,p) := \ppp_tv - \Delta v + (A\cdot\nabla)v + (v\cdot\nabla)B 
+ \nabla p.
$$
Then
$$
L(w, \wva^{\frac{m}{2}}p)(x,t) = \wva^{\frac{m}{2}}F + q(t)\wva w \quad 
\mbox{in $Q$}.
$$
Therefore Lemma 1 with $m=0$ yields 
\begin{align*}
& \int_Q (\vert \nabla w\vert^2 + s\va\vert \rrr w\vert^2
+ s^2\va^2\vert w\vert^2) \weight dxdt 
= \int_Q (\wva^m\vert \nabla v\vert^2 + s\wva^m\va\vert \rrr v\vert^2
+ s^2\wva^m\va^2\vert v\vert^2) \weight dxdt\\
\le &C\int_Q \vert \wva^{\frac{m}{2}}F + q\wva w\vert^2 \weight dxdt
+ C\int_{Q_{\ome}} (s\va\vert \rrr w\vert^2 + s^2\va^2\vert w\vert^2
+ s\va\vert \nabla w\vert^2) \weight dxdt \\
\le &C\int_Q \wva^m \vert F\vert^2 \weight dxdt
+ C\int_Q \wva^{m+2} \vert v\vert^2 \weight dxdt\\
+& C\int_{Q_{\ome}} (s\wva^{m+1}\vert \rrr v\vert^2 
+ s^2\wva^{m+2}\vert v\vert^2 + s\wva^{m+1}\vert \nabla v\vert^2) \weight dxdt
\end{align*}
for all large $s>0$.
Choosing $s>0$ large, in view of (2.8),
we can absorb the second term on the right-hand side
into the left-hand side, we obtain
\begin{align*}
& \int_Q (\wva^m\vert \nabla v\vert^2 + s\wva^{m+1}\vert \rrr v\vert^2
+ s^2\wva^{m+2}\vert v\vert^2) \weight dxdt\\
\le &C\int_Q \wva^m \vert F\vert^2 \weight dxdt
+ C\int_{Q_{\ome}} (s\wva^{m+1}\vert \rrr v\vert^2 
+ s^2\wva^{m+2}\vert v\vert^2 + s\wva^{m+1}\vert \nabla v\vert^2) \weight dxdt
\end{align*}
for all large $s>0$.
In terms of (2.8), we complete the proof of Lemma 1 with each $m\in \N$.
$\blacksquare$

\section*{Acknowledgment}
The first author was supported partly by NSF grant DMS 1312900.
The second author was supported by Grant-in-Aid for Scientific Research (S)
15H05740 and Grant-in-Aid (A) 20H00117 of 
Japan Society for the Promotion of Science, 
The National Natural Science Foundation of China
(no. 11771270, 91730303), and the RUDN University 
Strategic Academic Leadership Program.

\end{document}